\documentclass[MSNbibl,number,citesort,seceqn,dvips]{arxbj}
\usepackage{mathbh}
\usepackage{mathrsfs}

\aid{0}
\volume{20}
\issue{4}
\pubyear{2014}
\firstpage{1698}
\lastpage{1716}
\doi{10.3150/13-BEJ537} %kopijuoti is PTS

\makeatletter
\newcommand{\RMe}{\mathrm{e}}

\newcommand{\mrmd}{\,\mathrm{d}}
\newcommand{\mrmdd}{\mathrm{d}}

\newcommand{\rrvert}{\vert}
\newcommand{\llvert}{\vert}

\newcommand{\xrightarrow}[1]{\stackrel{#1}{\rightarrow}}

\newcommand{\iint}{\int\!\!\int}
\newcommand{\idotsint}{\int\cdots\int}

\newcommand{\ee}{\mathbb{E}}
\newcommand{\nn}{\mathbb{N}}
\newcommand{\rr}{\mathbb{R}}
\newcommand{\pp}{\mathbb{P}}

\def\MM{\mathcal M}

\newtheorem{theorem}{Theorem}[section]
\newtheorem{lem}[theorem]{Lemma}
\newremark{remarks}[theorem]{Remark}

\makeatother

\begin{document}
\begin{frontmatter}

\title{Convergence rate and concentration inequalities for Gibbs
sampling in high dimension}
\runtitle{Convergence rate and concentration inequalities}

\begin{aug}
%%%% inicialai - be tarpu
\author[1]{\inits{N.-Y.}\fnms{Neng-Yi} \snm{Wang}\corref{}\thanksref{1}\ead[label=e1]{wangnengyi@amss.ac.cn}} \and
\author[2]{\inits{L.}\fnms{Liming} \snm{Wu}\thanksref{2}\ead[label=e2]{Li-Ming.Wu@math.univ-bpclermont.fr}}
%%\runauthor{} %% auto
\address[1]{Institute of Applied Mathematics, Academy of Mathematics
and Systems Science, Chinese Academy of Sciences,
100190, Beijing, China. \printead{e1}}
\address[2]{Institute of Applied Mathematics, Academy of Mathematics
and Systems Science, Chinese Academy of Sciences,
100190, Beijing, China and Laboratoire de Math. CNRS-UMR 6620,
Universit\'e Blaise Pascal, 63177 Aubi\`ere, France. \printead{e2}}
\end{aug}

% HISTORY:
\received{\smonth{2} \syear{2013}}
\revised{\smonth{4} \syear{2013}}

% ABSTRACT
%
\begin{abstract}
The objective of this paper is to study the Gibbs sampling for
computing the mean of observable in very high dimension -- a powerful
Markov chain Monte Carlo method. Under the Dobrushin's uniqueness
condition, we establish some explicit and sharp estimate of the
exponential convergence rate and prove some Gaussian concentration
inequalities for the empirical mean.
\end{abstract}

% KEYWORDS
% visi is mazosios raides ir pagal abecele
%

\end{frontmatter}

%s1 #&#
\section{Introduction}\label{sec1}

Let $\mu$ be a Gibbs probability measure on $E^{N}$ with dimension
$N$ very big, that is,
\[
\mu\bigl(\mrmdd x^{1},\ldots,\mrmdd x^{N}\bigr)=\frac{\RMe ^{-V(x^{1},\ldots,x^{N})}}{\idotsint
_{E^{N}}\RMe^{-V(x^{1},\ldots,x^{N})}\pi(\mrmdd x^{1})\cdots
\pi(\mrmdd x^{N})}\pi
\bigl(\mrmdd x^{1}\bigr)\cdots\pi\bigl(\mrmdd x^{N}\bigr),
\]
where $\pi$ is some
$\sigma$-finite reference measure on $E$. Our purpose is to study
the Gibbs sampling -- a Markov chain Monte Carlo method (MCMC in
short) for approximating $\mu$. In fact, even for the simplest case
where $E=\{+, -\}$, as the denominator contains an exponential
number of terms and each of them may be very big or small for high
dimension, it is very difficult to model $\mu$.

Let $\mu_{i}(\cdot|x)$ $(x=(x^1,\ldots, x^N)\in E^{N})$ be the
regular conditional distribution of $x^{i}$ knowing $(x^{j},j\neq
i)$ under $\mu$; and $\bar{\mu}_{i}(\mrmdd y|x)= (\prod_{j\neq
i}\delta_{x^{j}}(\mrmdd y^{j}) )\otimes\mu_{i}(\mrmdd y^{i}|x)$
(product measure), where $\delta_\cdot$ is the Dirac measure at the point
$\cdot$. We see that
\[
\mu_i\bigl(\mrmdd x^i|x\bigr) =\frac{\RMe^{-V(x^{1},\ldots,x^{N})}}{\int
_{E}\RMe^{-V(x^{1},\ldots,x^{N})}\pi(\mrmdd x^{i})}\pi
\bigl(\mrmdd x^i\bigr),
\]
which is a one-dimensional measure, easy to be realized in practice.

The idea of the Gibbs sampling consists in approximating $\mu$ via
iterations of the one-dimensional conditional distributions
$\mu_{i}, i=1,\ldots, N$. It is described as follows. Given a
starting configuration $x_{0}=(x_{0}^{1},\ldots,x_{0}^{N})\in
E^{N}$, let $(X_{n}; {n\geq0})$ be a non-homogeneous Markov chain
defined on some probability space $(\Omega, \mathscr{F},
\mathbb{P}_{x_0})$, such that
\[
X_{0}=x_{0}\qquad \mbox{a.s.},
\]
and given
\[
X_{kN+i-1}=x=\bigl(x^{1},\ldots,x^{N}\bigr)\in
E^{N}\qquad (k\in\mathbb{N},1\leq i\leq N), %
\]
then $X_{kN+i}^{j}=x^{j}$ for $j\neq i$ and the conditional law of
$X_{kN+i}^{i}$ is $\mu_{i}(\mrmdd y^{i}|x)$. In other words, the transition
probability at step $kN+i$ is:
\begin{itemize}
\item[$(1)$]
$\mathbb{P}(X_{kN+i}\in \mrmdd y|X_{kN+i-1}=x)=\bar{\mu}_{i}(\mrmdd y|x)$.
\end{itemize}
Therefore,
\begin{itemize}
\item[$(2)$] $\mathbb{P}(X_{kN}\in \mrmdd y|X_{(k-1)N}=x)=(\bar{\mu
}_{1}\cdots\bar{\mu}_{N})(x,\mrmdd y)=:P(x,\mrmdd y)$.
\end{itemize}

Finally, the Gibbs sampling is the time-homogeneous Markov chain
$(Z_{k}=X_{kN},k=0,1,\ldots)$, whose transition probability is $P$.

This MCMC algorithm is known sometimes as \textit{Gibbs sampler} in the
literature (see Winkler~\cite{Winkler}, Chapters~5 and 6). It is
actively used in statistical physics, chemistry, biology and
throughout the Bayesian statistics (a sentence taken from
\cite{DKS09}). It was used by Zegarlinski \cite{Zeg92} as a tool
for proving the logarithmic Sobolev inequality for Gibbs measures,
see also the second named author \cite{Wu06} for a continuous time
MCMC.

Our purpose is two-fold:
\begin{itemize}
\item[$(1)$] the convergence rate of $P^{k}$ to $\mu$;

\item[$(2)$]
the concentration inequality for
$\mathbb{P}(\frac{1}{n}\sum_{k=1}^{n}f(Z_{k})-\mu(f)\geq t), t>0$.
\end{itemize}

Question (1) is a classic subject. Earlier works by Meyn and Tweedie
\cite{MT93} and Rosenthal \mbox{\cite{Rosenthal95,Rosenthal96}} are based
on the Harris ergodicity theorem (minorization condition together
with the drift condition in the non-compact case). Quantitative
estimates in the Harris ergodic theorem are obtained more recently
by Rosenthal \cite{Rosenthal02} and Hairer and Mattingly \cite{HM11}.
But as indicated by Diaconis, Khare and Saloff-Coste \cite
{DKS08,DKS09}, theoretical results obtained from the Harris theorem are
very far (even too far) from the convergence rate of numerical
simulations in high dimension (e.g., $N=100$). That is why Diaconis,
Khare and Saloff-Coste \cite{DKS08,DKS09} use new methods and tools
(orthogonal polynomials, stochastic monotonicity and coupling) for
obtaining sharp estimates of $\|\nu P^k -\mu\|_{\mathrm{TV}}$ (total
variation norm) for several special models in Bayesian statistics,
with $E^N$ replaced by $E\times\Theta$, a space of two different
components.

For the question (1), our tool will be the Dobrushin interdependence
coefficients (very natural and widely used in statistical physics),
instead of the minorization condition in the Harris theorem or the
special tools in \cite{DKS08,DKS09}. Our main idea consists in
constructing an appropriate coupling well adapted to the Dobrushin
interdependence coefficients, close to that of Marton
\cite{Marton04}.

To the second question, we will apply the recent theory on transport
inequalities (see Marton \cite{Marton96}, Ledoux \cite
{Ledoux99,Ledoux01}, Villani \cite{Vill03}, Gozlan and L\'eonard \cite{GL10}
and references therein), and our approach is inspired from Marton
\cite{Marton97,Marton04} and Djellout, Guillin and Wu \cite{DGW} for
dependent tensorization of transport inequalities.

See \cite{Winkler,DFG01,RR04} for Monte Carlo algorithms and
diverse applications, and \cite{JO} for concentration inequalities
of general MCMC under the positive curvature condition.

This paper is organized as follows. The main results are stated in
the next section, and we prove them in Section~\ref{sec3}.

%s2 #&#
\section{Main results}\label{sec2}
Throughout the paper, $E$ is a Polish space with the Borel
$\sigma$-field $\mathscr{B}$, and $d$ is a metric on $E$ such that
$d(x,y)$ is lower semi-continuous on $E^2$ (so $d$ does not
necessarily generate the topology of~$E$). On the product space we
consider the $L^1$-metric
\[
d_{L^{1}}(x,y):=\sum_{i=1}^{N}d
\bigl(x^{i},y^{i}\bigr),\qquad x,y\in E^{N}.
\]
If $d(x^i,y^i)=1_{x^i\ne y^i}$ is the discrete metric on $E$,
$d_{L^1}$ becomes the Hamming distance on $E^N$, a good metric
for concentration
 in
high dimension as shown by Marton \cite{Marton96,Marton97}.

%s2.1 #&#
\subsection{Dobrushin's interdependence coefficient}\label{sec2.1} Let $\MM_{1}(E)$
be the space
of probability measures on $E$ and
\[
\MM_{1}^{d}(E):= \biggl\{\nu\in\MM_{1}(E);
\int_{E}d(x_{0},x)\nu(\mrmdd x)< \infty\biggr\}
\]
$(x_{0}\in E
$ is some fixed point). Given $\nu_{1}, \nu_{2}\in\MM_{1}^{d}(E)$,
the $L^{1}$-Wasserstein distance between $\nu_{1}, \nu_{2}$ is given
by
%
%e2.1 #&#
%
\begin{equation}
W_{1,d}(\nu_{1},\nu_{2}):=\inf
_{\pi} \iint_{E\times
E}d(x,y)\pi(\mrmdd x,\mrmdd y),
\end{equation}
where the infimum is taken over all probability measures $\pi$ on
$E\times E$ such that its marginal distributions are, respectively,
$\nu_{1}$ and $\nu_{2}$ (\textit{coupling of $\nu_{1}$ and $\nu_{2}$,
say}). When $d(x,y)=\mathbh{1}_{x\neq y}$ (\textit{the discrete
metric}), it is well known that
\[
W_{1,
d}(\nu_{1},\nu_{2})=\sup
_{A\in
\mathscr{B}}\bigl|\nu_{1}(A)-\nu_{2}(A)\bigr|=
\frac{1}{2}\|\nu_{1}-\nu_{2}\| _{\mathrm{TV}}\qquad
\mbox{(total variation)}.
\]

Recall the Kantorovich--Rubinstein duality
relation \cite{Vill03}
\[
W_{1,d}(\nu_{1},\nu_{2})=\sup
_{\|f\|_{\mathrm{Lip}}\leq
1}\int_{E}fd(\nu_{1}-
\nu_{2}),\qquad \|f\|_{\mathrm{Lip}}:=\sup_{x\neq
y}
\frac{|f(x)-f(y)|}{d(x,y)}.
\]
Let $\mu_{i}(\mrmdd x^{i}|x)$ be the given
regular conditional distribution of $x^{i}$ knowing $(x^{j},j\neq
i)$.

\textit{Throughout the paper}, \textit{we assume that $\int_{E^N} d(y^i, x_0^i)
\mrmd \mu(y)<\infty$, $\mu_i(\cdot|x)\in\MM_1^{d}(E)$ for all
$i=1,\ldots,N$ and $x\in E^N$, where $x_0$ is some fixed point of
$E^N$, and $x\to\mu_i(\cdot|x)$ is Lipschitzian from
$(E^N,d_{L^1})$ to $(\MM_1^{d}(E), W_{1,d})$.}

Define the matrix of the $d$-Dobrushin interdependence
coefficients $C:=(c_{ij})_{i,j=1,\ldots,N}$
%
%e2.2 #&#
%
\begin{equation}
\label{coefficient} c_{ij}:=\sup_{x=y\ \mathrm{off}\ j}
\frac{W_{1,d}(\mu_{i}(\cdot|x),\mu
_{i}(\cdot|y))}{d(x^{j},y^{j})},\qquad i,j=1,\ldots,N.
\end{equation}
Obviously $c_{ii}=0$. Then the well-known Dobrushin uniqueness
condition (see \cite{Dobrushin68,Dobrushin70}) is read as
\[
\mbox{(H1)}\quad r:=\|C\|_{\infty}:=\max_{1\leq i\leq N}\sum
_{j=1}^{N}c_{ij}< 1
\]
or
\[
\mbox{(H2)}\quad r_1:=\|C\|_1:=\max_{1\leq j\leq N}\sum
_{i=1}^{N}c_{ij}< 1.
\]

By the triangular inequality
for the metric $W_{1,d}$,
%
%e2.3 #&#
%
\begin{equation}
\label{Dobrushin} W_{1,d}\bigl(\mu_{i}(\cdot|x),
\mu_{i}(\cdot|y)\bigr)\leq\sum_{j=1}^{N}c_{ij}d
\bigl(x^{j},y^{j}\bigr),\qquad 1\leq i\leq N.
\end{equation}

%s2.2 #&#
\subsection{Transport inequality and Bobkov--G\"{o}tze's criterion}\label{sec2.2}
When $\mu,\nu$ are probability measures, the Kullback information
(or relative entropy) of $\nu$ with respect to $\mu$ is defined as
%
%e2.4 #&#
%
\begin{equation}
H(\nu|\mu)=\cases{\displaystyle \int\log\frac{\mrmdd \nu}{\mrmdd \mu}\mrmd \nu, &\quad if $\nu\ll\mu$,
\vspace*{2pt}\cr
+\infty,
&\quad otherwise.}
\end{equation}

We say that the probability measure $\mu$ satisfies the
\textit{$L^{1}$-transport-entropy inequality} on $(E,d)$ with some
constant $C>0$, if
%
%e2.5 #&#
%
\begin{equation}
W_{1,d}(\mu,\nu)\leq\sqrt{2CH(\nu|\mu)},\qquad \nu\in\MM_1(E).
\end{equation}
To be short, we write $\mu\in T_{1}(C)$ for
this relation. This inequality, related to the phenomenon of measure
concentration, was introduced and studied by Marton \cite
{Marton96,Marton97}, developed subsequently by Talagrand \cite{Talagrand},
Bobkov and G\"{o}tze \cite{BG}, Djellout, Guillin and Wu \cite{DGW} and
amply explored by Ledoux \cite{Ledoux01,Ledoux99}, Villani
\cite{Vill03} and Gozlan-L\'eonard \cite{GL10}. Let us mention the
following Bobkov--G\"{o}tze's
criterion.
%
%le2.1 #&#
\begin{lem}[(\cite{BG})]\label{BG}
A probability measure $\mu$ satisfies the $L^{1}$-transport-entropy
inequality on $(E,d)$ with constant $C>0$, that is, $\mu\in
T_{1}(C)$, if and only if for any Lipschitzian function
$F\dvtx (E,d)\rightarrow\mathbb{R}$, $F$ is $\mu$-integrable and
\[
\int_{E}\RMe^{\lambda(F-\langle F\rangle_{\mu})}\mrmd \mu\leq\exp\biggl\{
\frac{\lambda^{2}}{2}C\|F\|_{\mathrm{Lip}}^{2}\biggr\}\qquad \forall\lambda\in
\mathbb{R},
\]
where $\langle F\rangle_{\mu}=\mu(F)=\int_{E}F \mrmd  \mu$. In that case,
\[
\mu\bigl(F-\langle F\rangle_{\mu}\geq t\bigr)\leq\exp\biggl\{-
\frac{t^{2}}{2C\|F\|_{\mathrm{Lip}}^{2}} \biggr\},\qquad t>0.
\]
\end{lem}

Another necessary and sufficient condition for $T_1(C)$ is the
Gaussian integrability of $\mu$, see Djellout, Guillin and Wu \cite{DGW}.
For further results and recent progresses see Gozlan and L\'eonard
\cite{GL07,GL10}.

%re2.2 #&#
\begin{remarks}\label{rem21}
Recall also that w.r.t. the discrete metric
$d(x,y)=1_{x\ne y}$, any probability measure $\mu$ satisfies
$T_1(C)$ with the sharp constant $C=1/4$ (the well known CKP
inequality).
\end{remarks}

%s2.3 #&#
\subsection{Main results}\label{sec2.3}
For any function $f\dvtx E^N\to\rr$, let
\[
\delta_i(f):=\sup_{x=y\ \mathrm{off}\ i} \frac{|f(x)-f(y)|}{d(x^i,y^i)}
\]
be the Lipschitzian coefficient w.r.t. the $i$th coordinate $x^i$.
It is easy to see that
\[
\|f\|_{\mathrm{Lip}(d_{L^1})}=\max_{1\le i\le N}\delta_i(f).
\]

%th2.3 #&#
\begin{theorem}[(Convergence rate)]\label{convergencerate}
Under the Dobrushin uniqueness condition \textup{(H1)}, we have:

\begin{enumerate}[(b)]
\item[(a)] For any Lipschitzian function $f$ on $E^N$ and two
initial distributions $\nu_1,\nu_2$
on $E^N$,
%
%e2.6 #&#
%
\begin{equation}
\label{a1} \bigl|\nu_1 P^kf -\nu_2
P^kf\bigr|\le r^k \max_{1\leq
i\leq N} \mathbb{E}d
\bigl(Z_{0}^{i}(1),Z_{0}^{i}(2)\bigr)
\sum_{i=1}^N \delta_i(f),\qquad k
\ge1,\quad
\end{equation}
where $(Z_{0}(1),Z_{0}(2))$ is a coupling
of $(\nu_{1},\nu_{2})$, that is, the law of $Z_0(j)$ is $\nu_j$ for
$j=1,2$.

\item[(b)] In particular for any initial distribution $\nu$ on $E^N$,
\[
W_{1,d_{L^{1}}}\bigl(\nu P^{k}, \mu\bigr)\leq Nr^{k}\max
_{1\leq i\leq N} \mathbb{E}d\bigl(Z_{0}^{i}(1),Z_{0}^{i}(2)
\bigr),
\]
where $(Z_{0}(1),Z_{0}(2))$ is a coupling of $(\nu,\mu)$.
\end{enumerate}
\end{theorem}

By part (b) above $\mu$ is the unique invariant measure of $P$ under
the Dobrushin uniqueness condition, and $P^k(x,\cdot)$ converges
exponentially rapidly to $\mu$ in the metric $W_{1,d_{L^1}}$,
showing theoretically why the numerical simulations by the Gibbs
sampling are very rapid.

%re2.4 #&#
\begin{remarks}
Let us compare Theorem \ref{convergencerate} with the
known results in \cite{MT93,Rosenthal95,Rosenthal96,DKS08,DKS09}
on the convergence rate of the Gibbs sampling.

At first the convergence rate in those known works is in the total
variation norm, not in the metric $W_{1,d_{L^1}}$. When $d$ is the
discrete metric, we have by part (b) of Theorem \ref{convergencerate}
\[
\tfrac12 \bigl\|\nu P^k-\mu\bigr\|_{\mathrm{TV}} \le W_{1,d_{L^{1}}}\bigl(\nu
P^{k}, \mu\bigr)\le Nr^k\qquad \forall\nu, \forall k\ge1.
\]

Next, let us explain once again why the minorization condition in the
Harris theorem does not yield accurate estimates in high dimension
(see Diaconis \textit{et al.} \cite{DKS08,DKS09} for similar
discussions based on concrete examples). Indeed assume that $E$ is
finite, then under reasonable assumption on $V(x^1,\ldots, x^N)$,
there are constant $c>0$ and a probability measure $\nu(\mrmdd y)$ such
that
\[
P(x,\mrmdd y)\ge \RMe^{-c N} \nu(\mrmdd y) %
\]
(i.e., almost the best minorization that one can obtain in the
dependent case). Hence by the Doeblin theorem (the ancestor of the
Harris theorem),
\[
\frac12 \sup_{x\in E^N}\bigl\|P^k(x,\cdot) -\mu
\bigr\|_{\mathrm{TV}} \le\bigl(1- \RMe^{-c
N}\bigr)^k. %
\]
So one requires at least an exponential number $\RMe^{c N}$ of steps
for the right-hand side becoming small. Our estimate of the
convergence rate is much better in high dimension, that is the good
point of Theorem \ref{convergencerate}.

The weak point of Theorem \ref{convergencerate} is that our result
depends on the Dobrushin uniqueness condition, even in low
dimension. If $N$ is small, the results in \cite
{MT93,Rosenthal95,Rosenthal96} are already good enough. Particularly
the estimates of
Diaconis, Khare and Saloff-Coste \cite{DKS08,DKS09} for the special
space of two different components in Bayesian statistics are sharp.

We should indicate that the Dobrushin uniqueness condition is quite
natural for the exponential convergence of $P^k$ to $\mu$ with the
rate $r$ independent of $N$ as in this theorem, since the Dobrushin
uniqueness condition is well known to be sharp for the phase
transition of mean field models \cite{Dobrushin68,Dobrushin70,DoSh85}.

Finally, our tool (Dobrushin's interdependence coefficients) is
completely different from those in the known works.
\end{remarks}

%re2.5 #&#
\begin{remarks}
As indicated by a referee, it would be very interesting to investigate
the convergence rate problem under the more flexible
Dobrushin--Shlosman analyticity condition (i.e., box version of
Dobrushin uniqueness condition, reference \cite{DoSh85}), but in that case
we feel that we should change the algorithm: instead of $\mu_i$, one
uses the conditional distribution $\mu_I(\mrmdd x_I|x)$ of $x_I$ knowing
$(x_j, j\notin I)$ where $I$ is a box containing $i$.
\end{remarks}

%re2.6 #&#
\begin{remarks}
A much more classical topic is Glauber dynamics associated with the
Gibbs measures in finite or infinite volume. We are content here to
mention only Zegarlinski \cite{Zeg92}, Martinelli and Olivieri
\cite{MaOli94}, and the lecture notes of Martinelli \cite{Martinelli99}
for a great number of references.
\end{remarks}

The convergence rate estimate above will be our starting point for
computing the mean $\mu(f)$, that is, to approximate $\mu(f)$ by the
empirical mean $\frac{1}{n}\sum_{k=1}^{n}f(Z_{k})$.

%th2.7 #&#
\begin{theorem}\label{concentrationinequalities} Assume
\[
r_1:=\|C\|_1=\max_{1\le j\le N} \sum
_{i=1}^N c_{ij}<\frac12 %
\]
and for some constant $C_{1}>0$,
\[
\mbox{\textup{(H3)}}\quad \forall i=1,\ldots, N, \forall x\in E^N\qquad \mu_{i}(
\cdot|x) \in T_{1}(C_{1}).
\]
(Recall that $C_1=1/4$ for the discrete metric
$d(x,y)=1_{x\ne y}$.)
Then for any Lipschitzian function $f$ on
$E^{N}$ with $\|f\|_{\mathrm{Lip}(d_{L^1})}\leq\alpha$, we have:
\begin{enumerate}[(b)]
\item[(a)]
%
%e2.7 #&#
%
\begin{equation}
\label{thm2a}\mathbb{P}_{x} \Biggl(\frac{1}{n}\sum
_{k=1}^{n}f(Z_{k})-\frac{1}{n}\sum
_{k=1}^{n}P^{k}f(x)\geq t \Biggr)
\leq\exp\biggl\{-\frac{t^{2}(1-2r_{1})^{2}n}{2C_{1}\alpha^{2}
N} \biggr\}\qquad \forall t>0,n\geq1;
\end{equation}
\item[(b)] furthermore if \textup{(H1)} holds,
%
%e2.8 #&#
%
\begin{eqnarray}
\label{thm2b}
&&
\mathbb{P}_{x} \Biggl(\frac{1}{n}\sum
_{k=1}^{n}f(Z_{k})-\mu(f)\geq
\frac{ M}{n}+t \Biggr)\nonumber\\[-8pt]\\[-8pt]
&&\quad\leq\exp\biggl\{-\frac
{t^{2}(1-2r_{1})^{2}n}{2C_{1}\alpha^{2}
N} \biggr\}\qquad \forall t>0,n\geq1,\nonumber
\end{eqnarray}
where
\[
M=\frac{r}{1-r}\max_{1\leq i\leq
N}\int_{E^{N}}d
\bigl(x^{i},y^{i}\bigr)\mrmd \mu(y)\cdot\sum
_{i=1}^N \delta_i(f).
\]
\end{enumerate}
\end{theorem}
In conclusion under the conditions of this theorem, when $n\gg N$,
the empirical means $\frac{1}{n}\sum_{k=1}^{n}f(Z_{k})$ will
approximate to $\mu(f)$ exponentially rapidly in probability with
the speed $n/N$, with the bias not greater than $M/n$. The speed
$n/N$ is the correct one, as will be shown in the remark below.

We do not know whether the concentration inequality with the speed
$n/N$ still holds under the more natural Dobrushin's uniqueness
condition $r_1=\|C\|_1<1$. We know only that $r_1<1$ does not imply
that $P$ is contracting in the metric $W_{1,d_{L^1}}$, see the
example in Remark \ref{remar3.3}.

%re2.8 #&#
\begin{remarks}
Consider $f(x)=\frac1N \sum_{i=1}^N g(x^i)$ where $g\dvtx E\to \rr$ is
$d$-Lipschitzian with $\|g\|_{\mathrm{Lip}}=\alpha$ (the observable of this type
is often used in statistical mechanics). Since
$\|f\|_{\mathrm{Lip}(d_{L^1})}=\frac{\alpha}{N}$, the inequality (\ref{thm2a})
implies for all $t>0,n\geq1$,
%
%e2.9 #&#
%
\begin{equation}
\label{thm2d} \mathbb{P}_{x} \Biggl(\frac{1}{nN}\sum
_{k=1}^{n}\sum_{i=1}^Ng
\bigl(Z_{k}^i\bigr)- \ee_x \frac{1}{nN}
\sum_{k=1}^{n}\sum
_{i=1}^Ng\bigl(Z_{k}^i\bigr)
\geq t \Biggr)\leq\exp\biggl\{-\frac
{t^{2}(1-2r_{1})^{2}nN}{2C_{1}\alpha^{2} } \biggr\},
\end{equation}
which is of speed $nN$.

Let us show that the concentration inequality (\ref{thm2a}) is
sharp. In fact in the free case, that is, $\mu_i(\mrmdd x_i|x)=\beta(\mrmdd x_i)$
does not depend upon $(x_j)_{j\ne i}$ and $i$, and $\mu$ is the
product measure $\beta^{\otimes N}$. In this case $P(x,\mrmdd y)=\mu(\mrmdd y)$,
in other words $(Z_k)_{k\ge1}$ is a sequence of independent and
identically distributed (i.i.d. in short) random variables valued in
$E^N$, of common law $\mu$. Since $r_1=0$ in the free case, the
concentration inequality (\ref{thm2d}) is equivalent to the
transport inequality \textup{(H3)} for $\mu_i=\beta$, by Gozlan-L\'eonard
\cite{GL07}. That shows also the speed $n/N$ in Theorem
\ref{concentrationinequalities} is the correct one.
\end{remarks}

%re2.9 #&#
\begin{remarks}
We explain now why we do not apply directly the nice concentration
results of Joulin and Ollivier \cite{JO} for general MCMC. In fact
under the condition that $r_1<1/2$, we can prove that
\[
W_{1,d_{L^1}}\bigl(P(x,\cdot), P(y,\cdot)\bigr)\le\frac{r_1}{1-r_1}
d_{L^1}(x,y)
\]
(by Lemma \ref{1-norm}). In other words the Ricci curvature in
\cite{JO} is bounded from below by
\[
\kappa:=1-\sup_{x\ne y}\frac{W_{1,d_{L^1}}(P(x,\cdot),
P(y,\cdot))}{d_{L^1}(x,y) } \ge1-
\frac
{r_1}{1-r_1}=\frac{1-2r_1}{1-r_1}.
\]
Unfortunately we cannot show
that the Ricci curvature $\kappa$ is positive in the case where
$r_1\in[1/2,1)$.

If $(E,d)$ is unbounded, the results of \cite{JO}, Theorems 4 and 5,
do not apply here, because their granularity constant
\[
\sigma_\infty= \frac12 \sup_{x\in E^N} \operatorname{Diam}\bigl(\operatorname{supp}
\bigl(P(x,\cdot)\bigr)\bigr) %
\]
explodes.\vadjust{\goodbreak}

Assume now that $(E,d)$ is bounded.
If we apply the results (\cite{JO}, Theorems 4 and 5) and their notations,
their coarse diffusion constant
\[
\sigma(x)^2=\frac12 \iint d_{L^1}(y,z)^2
P(x,\mrmdd y) P(x,\mrmdd z) %
\]
is of order $N^2$; and their local dimension
\[
n_x = \inf_{f: \|f\|_{\mathrm{Lip}(d_{L^1})}=1}
\frac{\sigma(x)^2}{\operatorname{Var}_{P(x,\cdot)}(f)} %
\]
is of order $N$ (by Lemma \ref{transine} below), and their
granularity constant $\sigma_\infty$ is of order
$N$. Setting
\[
V^2 = \frac{1}{\kappa n } \sup_{x\in E^N}
\frac{\sigma(x)^2}{n_x
\kappa};\qquad r_{\max}=\frac{4V^2\kappa n}{3\sigma_\infty}.
\]
Theorem 4 in \cite{JO} says that if $\|f\|_{\mathrm{Lip}(d_{L^1})}=1$,
\[
\mathbb{P}_{x} \Biggl(\Biggl\llvert\frac{1}{n}\sum
_{k=1}^{n}f(Z_{k})-\frac
{1}{n}\sum
_{k=1}^{n}P^{k}f(x)\Biggr\rrvert
\geq t \Biggr)\leq\cases{\displaystyle 2 \exp\biggl(-\frac{t^{2}}{16
V^2 } \biggr), &\quad $t
\in(0,r_{\max})$,
\vspace*{2pt}\cr
\displaystyle 2 \exp\biggl(-\frac{\kappa n t}{12 \sigma_\infty} \biggr), &\quad $t\ge
r_{\max}$,} %
\]
for all $n\ge1$ and $t>0$. So for small deviation $t$, their result
yields the same order Gaussian concentration inequality, but for
large deviation $t$, their estimate is only exponential, not
Gaussian as one may expect in this bounded case. In
\cite{JO}, Theorem 5, they get a same type Gaussian-exponential concentration
inequality with $V^2, r_{\max}$ depending upon the starting point
$x$.

Anyway the key lemmas in this paper are necessary for applying the
results of \cite{JO} to this particular model.
\end{remarks}

%re2.10 #&#
\begin{remarks}
For the Gibbs measure $\mu$ on $\rr^N$, Marton \cite{Marton04}
established the Talagrand transport inequality $T_2$ on $\rr^N$
equipped with the Euclidean metric, under the Dobrushin--Shlosman
analyticity type condition. The second named author \cite{Wu06} proved
$T_1(C)$ for $\mu$ on $E^N$ equipped with the metric $d_{L^1}$, under
\textup{(H1)}. But those transport inequalities are for the equilibrium
distribution~$\mu$, not for the Gibbs sampling which is a Markov chain
with $\mu$ as invariant measure. However our coupling is very close to
that of K. Marton.
\end{remarks}

%re2.11 #&#
\begin{remarks}
For $\phi$-mixing sequence of dependent random variables, Rio
\cite{Rio00} and Samson \cite{Samson00} established accurate
concentration inequalities, see also Djellout, Guillin and Wu
\cite{DGW} and the recent works by Paulin \cite{Paulin12} and
Wintenberger \cite{Winten12} for generalizations and improvements. In
the Markov chain case $\phi$-mixing means the Doeblin uniform
ergodicity. If one applies the results in \cite{Rio00,Samson00} to the
Gibbs sampling, one obtains the concentration inequalities with the
speed $n/(NS^2)$, where
\[
S=\sum_{k=0}^\infty\sup_{x\in E^N}
\bigl\|\nu P^k -\mu\bigr\|_{\mathrm{TV}}.
\]
When
\textup{(H1)} holds with the discrete metric $d$, $S$ is actually finite but
it is of order $N$ by Theorem~\ref{convergencerate} (and its
remarks). The concentration inequalities so obtained from
\cite{Rio00,Samson00} are of speed $n/N^3$, very far from the
correct speed $n/N$.
\end{remarks}

%re2.12 #&#
\begin{remarks}
When $f$ depends on a very small number of variables, since
$\|f\|_{\mathrm{Lip}(d_{L^1})}=\max_i \delta_i(f)$ does not reflect the nature
of such observable, one can imagine that our concentration inequalities
do not yield the correct speed. In fact in the free case and for
$f(x)=g(x^1)$, the correct speed must be $n$, not $n/N$. For this type
of observable, one may use the metric $d_{L^2}$ which reflects much
better the number of variables in such observable. The ideas in
Marton \cite{Marton03,Marton04} should be helpful. That will be another history.
\end{remarks}

%s3 #&#
\section{Proofs of the main results}\label{sec3}

%s3.1 #&#
\subsection{The construction of the coupling}\label{sec3.1}
Given any two initial distributions $\nu_{1}$ and $\nu_{2}$ on $E^{N}$,
we begin by constructing our coupled non-homogeneous Markov chain
$(X_{i},Y_{i})_{i\geq0}$, which is quite close to the coupling by
Marton \cite{Marton04}.

Let
$(X_{0},Y_{0})$ be a coupling of $(\nu_{1},\nu_{2})$.
And
given
\[
(X_{kN+i-1},Y_{kN+i-1})=(x,y)\in E^{N}\times
E^{N},\qquad k\in\mathbb{N},1\leq i\leq N,
\]
then
\[
X_{kN+i}^{j}=x^{j},\qquad Y_{kN+i}^{j}=y^{j},\qquad
j\neq i,
\]
and
\[
\mathbb{P}\bigl(\bigl(X_{kN+i}^{i}, Y_{kN+i}^{i}
\bigr)\in\cdot|(X_{kn+i-1},Y_{kn+i-1})=(x,y)\bigr)=\pi(\cdot|x,y),
\]
where $\pi(\cdot|x,y)$ is an optimal coupling of $\mu_{i}(\cdot|x)$
and $\mu_{i}(\cdot|y)$ such that
\[
\iint_{E^{2}}d(\tilde{x},\tilde{y})\pi(\mrmdd \tilde{x},\mrmdd \tilde
{y}|x,y)=W_{1,d}\bigl(\mu_{i}(\cdot|x),\mu_{i}(
\cdot|y)\bigr).
\]
Define the partial order on $\mathbb{R}^{N}$ by $a\leq b$ if and
only if $a^{i}\leq b^{i}, i=1,\ldots,N$. Then, by
(\ref{Dobrushin}), we have for $\forall k\in\nn, 1\leq i\leq N$,
\[
\pmatrix{ \mathbb{E}\bigl[d\bigl(X_{kN+i}^{1},Y_{kN+i}^{1}
\bigr)|X_{kN+i-1}, Y_{kN+i-1}\bigr]
\cr
\vdots
\cr
\vdots
\cr
\mathbb{E}
\bigl[d\bigl(X_{kN+i}^{N},Y_{kN+i}^{N}
\bigr)|X_{kN+i-1}, Y_{kN+i-1}\bigr]} \leq B_{i} \pmatrix{ d
\bigl(X_{kN+i-1}^{1},Y_{kN+i-1}^{1}\bigr)
\cr
\vdots
\cr
\vdots
\cr
d\bigl(X_{kN+i-1}^{N},Y_{kN+i-1}^{N}
\bigr)},
\]
where
\[
B_{i}=\pmatrix{ 1 & & & & &
\cr
& \ddots& & & &
\cr
& & 1 & & &
\cr
c_{i1} & c_{i2} & \cdots& \cdots& \cdots& c_{iN}
\cr
& & & 1 & &
\cr
& & & & \ddots&
\cr
& & & & & 1} \qquad\mbox{(the blank in the matrix means 0)}.
\]

Therefore by iterations, we have
%
%e3.1 #&#
%
\begin{equation}
\label{Matrix} \pmatrix{ \mathbb{E}d\bigl(X_{N}^{1},Y_{N}^{1}
\bigr)
\cr
\vdots
\cr
\vdots
\cr
\mathbb{E}d\bigl(X_{N}^{N},Y_{N}^{N}
\bigr)} \leq B_{N} B_{N-1}\cdots B_{1} \pmatrix{
\mathbb{E}d\bigl(X_{0}^{1},Y_{0}^{1}
\bigr)
\cr
\vdots
\cr
\vdots
\cr
\mathbb{E}d\bigl(X_{0}^{N},Y_{0}^{N}
\bigr)}.
\end{equation}

Let
%
%e3.2 #&#
%
\begin{equation}
\label{Q} Q:=B_{N}B_{N-1}\cdots B_{1}.
\end{equation}
Then we have the following lemma.
%
%le3.1 #&#
\begin{lem}\label{infty-norm}
Under \textup{(H1)}, $ \|Q\|_{\infty}:=\max_{1\leq i\leq N}\sum
_{j=1}^{N}Q_{ij}\leq r$.
\end{lem}
\begin{pf}
We use the probabilistic method. Under \textup{(H1)} we can construct Markov
chain $\{ \xi_{0},\ldots,\xi_{N}\}$, taking values in
$\{1,\ldots,N\}\sqcup\Delta$ where $\Delta$ is an extra point
representing the cemetery, and write as follows:
\[
\xi_{0}\xrightarrow{B_{N}}\xi_{1}
\xrightarrow{B_{N-1}}\cdots\xrightarrow{B_{1}}
\xi_{N},
\]
where the transition matrix
from $\xi_{i}$ to $\xi_{i+1}$ is $B_{N-i}$, more precisely
for $\forall
k=1,\ldots, N$,
\begin{eqnarray*}
\mathbb{P}\bigl(\xi_{k}=j|\xi_{k-1}=N-(k-1)
\bigr)&=&c_{N-(k-1),j},
\\
\mathbb{P}\bigl(\xi_{k}=\Delta|\xi_{k-1}=N-(k-1)\bigr)&=&1-
\sum_{j=1}^{N}c_{N-(k-1),j},
\\
\mathbb{P}(\xi_{k}=j|\xi_{k-1}=i)&=&\delta_{ij},\qquad i \neq N-(k-1),
\\
\mathbb{P}(\xi_{k}=\Delta|\xi_{k-1}=\Delta)&=&1.
\end{eqnarray*}
Here $\delta_{ij}=1$ if $i=j$ and $0$ otherwise (Kronecker's symbol).
Then
\[
\|Q\|_{\infty} =\max_{1\leq i\leq
N}\mathbb{P}(\xi_{N}
\neq\Delta|\xi_0=i).
\]
For any $i=1,\ldots,N$, when $\xi_{0}=i$, we have
$\xi_{0}=\xi_{1}=\cdots=\xi_{N-i}=i$. Therefore,
\[
\mathbb{P}(\xi_{N-i+1}\neq\Delta|\xi_0=i)=\sum
_{j=1}^{N}c_{ij}\leq r\vspace*{-3pt}
\]
and thus
\[
\mathbb{P}(\xi_{N}\neq\Delta|\xi_0=i)\leq\mathbb{P}(
\xi_{N-1}\neq\Delta|\xi_0=i)\leq\cdots\leq\mathbb{P}(
\xi_{N-i+1}\neq\Delta|\xi_0=i)\leq r.
\]
So $\|Q\|_{\infty} \leq
r$.\vspace*{-3pt}
\end{pf}

%s3.2 #&#
\subsection{Proof of Theorem
\texorpdfstring{\protect\ref{convergencerate}}{2.3}}\label{sec3.2}\vspace*{-3pt}

By (\ref{Matrix})
above, Markov property and iterations,
%
%e3.3 #&#
%
\begin{equation}
\label{matrixineq} \pmatrix{ \mathbb{E}d\bigl(X_{kN}^{1},Y_{kN}^{1}
\bigr)
\cr
\vdots
\cr
\vdots
\cr
\mathbb{E}d\bigl(X_{kN}^{N},Y_{kN}^{N}
\bigr)} \leq Q^{k} \pmatrix{ \mathbb{E}d\bigl(X_{0}^{1},Y_{0}^{1}
\bigr)
\cr
\vdots
\cr
\vdots
\cr
\mathbb{E}d\bigl(X_{0}^{N},Y_{0}^{N}
\bigr)}.
\end{equation}

Let $Z_{k}(1)=X_{kN}, Z_{k}(2)=Y_{kN}, k\geq0$, then by Lemma
\ref{infty-norm}
%
%e3.4 #&#
%
\begin{equation}
\label{maxine} \max_{1\leq i\leq
N}\mathbb{E}d\bigl(Z_{k}^{i}(1),
Z_{k}^{i}(2)\bigr)\leq r^{k} \max
_{1\leq
i\leq N}\mathbb{E}d\bigl(Z_{0}^{i}(1),
Z_{0}^{i}(2)\bigr).
\end{equation}
Now
the results of this theorem follow quite easily from this
inequality. In fact,

(a) For any Lipschitzian function $f\dvtx E^N\to\rr$,
\begin{eqnarray*}
\bigl|\nu_1 P^k f-\nu_2P^kf\bigr|
&=&\bigl|\ee f\bigl(Z_k(1)\bigr)-\ee f\bigl(Z_k(2)\bigr)\bigr|
\\
&\le&\sum_{i=1}^N \delta_i(f)
\ee d\bigl(Z_k^i(1), Z_k^i(2)
\bigr)
\\[-3pt]
&\le& r^k \max_{1\le i\le N}\ee d\bigl(Z_0^i(1),
Z_0^i(2)\bigr) \sum_{i=1}^N
\delta_i(f),
\end{eqnarray*}
where the last inequality follows by (\ref{maxine}). That is (\ref{a1}).

(b) Now for $\nu_1=\nu$, $\nu_2=\mu$, as $\mu P=\mu$, we have
\begin{eqnarray*}
W_{1,d_{L^{1}}}\bigl(\nu P^{k}, \mu\bigr)&=&W_{1,d_{L^{1}}}\bigl(\nu
P^{k}, \mu P^k\bigr) \leq \mathbb{E}\sum
_{i=1}^{N}d\bigl(Z_{k}^{i}(1),
Z_{k}^{i}(2)\bigr)
\\
&\leq& N \max_{1\leq i\leq N}\mathbb{E}d\bigl(Z_{k}^{i}(1),
Z_{k}^{i}(2)\bigr)
\\
&\leq& N r^{k}\max_{1\leq i\leq N}\mathbb{E}d
\bigl(Z_{0}^{i}(1), Z_{0}^{i}(2)
\bigr),
\end{eqnarray*}
the desired result.\vadjust{\goodbreak}

%s3.3 #&#
\subsection{Proof of Theorem \texorpdfstring{\protect\ref{concentrationinequalities}}{2.7}}\label{sec3.3}
We begin with
%
%le3.2 #&#
\begin{lem}\label{1-norm} If $r_1:=\|C\|_1<1$ (i.e., \textup{(H2)}), then for
the matrix $Q$ given in (\ref{Q}),
%
%e3.5 #&#
%
\begin{equation}
\label{Q1}\|Q\|_{1}:=\max_{1\leq j\leq N}\sum
_{k=1}^{N}Q_{kj}\leq\frac{r_{1}}{1-r_{1}}.
\end{equation}
In particular
%
%e3.6 #&#
%
\begin{equation}
\label{Q2} W_{1,d_{L^1}}\bigl(P(x,\cdot), P(y,\cdot)\bigr)\le
\frac
{r_1}{1-r_1} d_{L^1}(x,y)\qquad \forall x,y\in E^N.
\end{equation}
\end{lem}

\begin{pf} The last conclusion (\ref{Q2}) follows from (\ref{Q1}) and
(\ref{Matrix}). We show now (\ref{Q1}).

By the definition of $Q=B_N\cdots B_1$, it is not difficult to
verify for $1\leq k\leq N$,
%
%e3.7 #&#
%
\begin{equation}
Q_{kj}=\cases{0, \qquad \mbox{if $j=1$},
\vspace*{2pt}\cr
\displaystyle \sum_{h=1}^{j-1}
\Biggl(\sum_{l=1}^{k-1} \sum
_{k>i_{l}>\cdots>i_{2}>i_{1}=h}c_{k,i_{l}}c_{i_{l},i_{l-1}}\cdots
c_{i_{2},i_{1}=h}c_{h,j}+c_{h,j}1_{h=k} \Biggr), \vspace*{2pt}\cr
\qquad\hspace*{9.85pt}
\mbox{if $2\leq j\leq N$}.}
\end{equation}
Here we make the convention $\sum_{\varnothing}\cdot=0$. This can be
obtained again by the Markov chain $(\xi_0, \xi_1,\ldots, \xi_N)$
valued in $\{1,\ldots, N\}\cup\{\Delta\}$ constructed in Lemma
\ref{infty-norm}. Since $\pp(\xi_N=1|\xi_{N-1} =i)=0$ for all $i$,
$Q_{k1}=\pp(\xi_N=1|\xi_0=k)=0$: that is the first line in the
expression of $Q$. Now for $j=2,\ldots, N$, as
\[
Q_{kj}=\pp(\xi_N=j|\xi_0=k)=\pp(
\xi_N=j|\xi_{N-k}=k) %
\]
and if $\xi_{N-k}=k $ and $\xi_{N-k+1}>k$, then
$\xi_{N-k+1}=\cdots=\xi_N$ and so
\[
\pp(\xi_N=j|\xi_{N-k}=k)= c_{kj}1_{k<j}
+\sum_{h=1}^{j-1} c_{kh} \pp(
\xi_N=j| \xi_{N-k+1}=h). %
\]
This implies the expression of $Q$ above by induction.

Thus for $2\leq j\leq N$,
\begin{eqnarray*}
\sum_{k=1}^{N}Q_{kj}&=&
\sum_{k=1}^{N}\sum
_{h=1}^{j-1} \Biggl(\sum_{l=1}^{k-1}
\sum_{k>i_{l}>\cdots>i_{2}>i_{1}=h}c_{k,i_{l}}c_{i_{l},i_{l-1}}\cdots
c_{i_{2},i_{1}=h}c_{h,j}+c_{h,j}1_{h=k} \Biggr)
\\
&=&\sum_{k=1}^{N}\sum
_{h=1}^{j-1}c_{hj}1_{h=k}+ \sum
_{k=1}^{N}\sum_{h=1}^{j-1}
\sum_{l=1}^{k-1}\sum
_{k>i_{l}>\cdots
>i_{2}>i_{1}=h}c_{k,i_{l}}c_{i_{l},i_{l-1}}\cdots
c_{i_{2},i_{1}=h}c_{h,j}
\\
&=&\sum_{h=1}^{j-1}c_{hj}+\sum
_{h=1}^{j-1}\sum
_{k=1}^{N}\sum_{l=1}^{k-1}
\sum_{k>i_{l}
>\cdots>i_{2}>i_{1}=h}c_{k,i_{l}}c_{i_{l},i_{l-1}}\cdots
c_{i_{2},i_{1}=h}c_{h,j}
\\
&\leq& r_{1}+\sum_{l=1}^{N-1}\sum
_{h=1}^{j-1}\sum
_{k=l+1}^{N}\sum_{k>i_{l}
>\cdots>i_{2}>i_{1}=h}c_{k,i_{l}}c_{i_{l},i_{l-1}}
\cdots c_{i_{2},i_{1}=h}c_{h,j}
\\
&\leq& r_{1}+r_{1}^{2}+\cdots+r_{1}^{N},
\end{eqnarray*}
where the last inequality holds because for fixed $l\dvt 1\leq l\leq
N-1$ and $h\dvt1\le h\le j-1$,
\[
\sum_{k=1}^{N}\sum
_{k>i_{l}
>\cdots>i_{2}>i_{1}=h}c_{k,i_{l}}c_{i_{l},i_{l-1}}\cdots
c_{i_{2},i_{1}=h}\leq\sum_{k=1}^N
\bigl(C^l\bigr)_{kh}\le\bigl\|C^l\bigr\|_1\le
r_{1}^{l}.
\]
So the proof of (\ref{Q1}) is completed.
\end{pf}

%re3.3 #&#
\begin{remarks}\label{remar3.3}
Let $\mu$ be the Gaussian distribution on $\rr^2$ with mean $0$ and the
covariance matrix $ \bigl({1 \atop r}\enskip{r \atop1} \bigr)$ where $r\in(0,1)$.
We have $c_{12}=c_{21}=r<1$ (i.e., \textup{(H1)} and \textup{(H2)} both hold); and under
$P((x_1,x_2), (\mrmdd y_1,\mrmdd y_2))$, $z_1=y_1- rx_2$ and $z_2=y_2-r y_1$ are
i.i.d. Gaussian random variables with mean $0$ and variance $1-r^2$.
Hence,
\[
Q=B_2 B_1 = \pmatrix{1 & 0
\cr
r &0} \pmatrix{0 & r
\cr
0
& 1} = \pmatrix{0 & r
\cr
0 & r}
\]
and since $y_1=rx_2+z_1$, $y_2=z_2 + rz_1 +r^2 x_2$,
\[
W_{1,d_{L^1}} \bigl[P\bigl((x_1,x_2), \cdot\bigr), P
\bigl(\bigl(x'_1,x'_2\bigr),
\cdot\bigr) \bigr]= \bigl(r+r^2\bigr) \bigl|x_2-x'_2\bigr|.
\]
Thus, $\|Q\|_1 =2r$ and the Ricci curvature $\kappa$ is positive if
and only if $r+r^2<1$. In other words, though we have missed many
terms in the proof above, the estimate of $\|Q\|_1$ cannot be
qualitatively improved.
\end{remarks}

%le3.4 #&#
\begin{lem}\label{transine}
Assume \textup{(H2)} and \textup{(H3)}, then
\[
P(x_0, \cdot)\in T_{1} \biggl(\frac{NC_{1}}{(1-r_{1})^{2}}
\biggr)\qquad \forall x_0=\bigl(x_0^{1},\ldots,x_0^{N}\bigr)\in E^{N}.
\]
\end{lem}

\begin{pf} The proof is similar to the one used by Djellout, Guillin and Wu
\cite{DGW}, Theorem 2.5.
First for simplicity denote $P(x_0,\cdot)$ by $P$ and note that for
$1\leq i\leq N$,
\begin{eqnarray*}
X_{N}^{1}&=&X_{1}^{1},\ldots,X_{N}^{i}=X_{i}^{i},
\\
P\bigl(X_{N}^{i}\in\cdot|X_{N}^{1},\ldots,X_{N}^{i-1}\bigr)&=&\mu_{i}\bigl(
\cdot|X_{N}^{1},\ldots,X_{N}^{i-1},x_0^{i+1},\ldots,x_0^{N}\bigr)
\end{eqnarray*}\eject\noindent
and thus
\begin{eqnarray*}
&&
W_{1,d}
\bigl(P\bigl(X_{N}^{i}\in
\cdot|X_{N}^{1}=x^{1},\ldots,X_{N}^{i-1}=x^{i-1}
\bigr),P\bigl(X_{N}^{i}\in\cdot|X_{N}^{1}=y^{1},\ldots,X_{N}^{i-1}=y^{i-1}\bigr)\bigr)\\
&&\quad
\leq\sum_{j=1}^{i-1}c_{ij}d
\bigl(x^{j},y^{j}\bigr).
\end{eqnarray*}
For any
probability measure $Q$ on $E^{N}$ such that $H(Q|P)<\infty$, let
$Q_{i}(\cdot|x^{[1,i-1]})$ be the regular conditional law of $x^{i}$
knowing $x^{[1,i-1]}$, where $i\geq
2,x^{[1,i-1]}=(x^{1},\ldots,x^{i-1})$, and
$Q_{i}(\cdot|x^{[1,i-1]})$ the law of $x^{1}$ for $i=1$, all under
law $Q$. Define $P_{i}(\cdot|x^{[1,i-1]})$ similarly but under $P$.
We shall use the Kullback information between conditional
distributions,
\[
H_{i}\bigl(y^{[1,i-1]}\bigr)=H\bigl(Q_{i}\bigl(
\cdot|y^{[1,i-1]}\bigr)|P_{i}\bigl(\cdot|y^{[1,i-1]}\bigr)
\bigr)
\]
and exploit the following important identity:
\[
H(Q|P)=\sum_{i=1}^{N}\int
_{E^{N}}H_{i}\bigl(y^{[1,i-1]}\bigr)\mrmd Q(y).
\]
The
key is to construct an appropriate coupling of $Q$ and $P$, that is,
two random sequences $Y^{[1,N]}$ and $X^{[1,N]}$ taking values on
$E^{N}$ distributed according to $Q$ and $P$, respectively, on some
probability space $(\Omega,\mathscr{F},\mathbb{P})$. We define a
joint distribution $\mathscr{L}(Y^{[1,N]},X^{[1,N]})$ by induction
as follows (the Marton coupling).

At first the law of $(Y^1,X^1)$ is the optimal coupling of $Q(x^1\in
\cdot)$ and $P(x^1\in\cdot)$ $(=\mu_1(\cdot|x_0))$.
Assume that for some $i,2\leq i\leq
N,(Y^{[1,i-1]},X^{[1,i-1]})= (y^{[1,i-1]},x^{[1,i-1]})$ is given.
Then the joint conditional distribution
$\mathscr{L}(Y^{i},X^{i}|Y^{[1,i-1]}=y^{[1,i-1]},X^{[1,i-1]}=x^{[1,i-1]})$
is the optimal coupling of $Q_i(\cdot|y^{[1,i-1]})$ and
$P_i(\cdot|x^{[1,i-1]})$, that is,
\[
\mathbb{E}\bigl(d\bigl(Y^{i},X^{i}\bigr
)|Y^{[1,i-1]}=y^{[1,i-1]},X^{[1,i-1]}=x^{[1,i-1]}
\bigr) = W_{1,d}\bigl(Q_{i}\bigl(\cdot|y^{[1,i-1]}
\bigr),P_{i}\bigl(\cdot|x^{[1,i-1]}\bigr)\bigr).
\]
Obviously, $Y^{[1,N]},X^{[1,N]}$ are of law $Q,P$, respectively. By
the triangle inequality for the $W_{1,d}$ distance,
\begin{eqnarray*}
&&\mathbb{E}\bigl(d\bigl(Y^{i},X^{i}
\bigr)|Y^{[1,i-1]}=y^{[1,i-1]},X^{[1,i-1]}=x^{[1,i-1]}\bigr)
\\
&&\quad\leq W_{1,d}\bigl(Q_{i}\bigl(\cdot|y^{[1,i-1]}
\bigr),P_{i}\bigl(\cdot|y^{[1,i-1]}\bigr)\bigr)+ W_{1,d}
\bigl(P_{i}\bigl(\cdot|y^{[1,i-1]}\bigr),P_{i}\bigl(
\cdot|x^{[1,i-1]}\bigr)\bigr)
\\
&&\quad\leq\sqrt{2C_{1}H_{i}\bigl(y^{[1,i-1]}\bigr)}+\sum
_{j=1}^{i-1}c_{ij}d
\bigl(x^{j},y^{j}\bigr).
\end{eqnarray*}
By recurrence on $i$, this entails that
$\mathbb{E}d(Y^{i},X^{i})<\infty$ for all $i=1,\ldots,N$. Taking the
average with respect to $\mathscr{L}(Y^{[1,i-1]},X^{[1,i-1]})$,
summing on $i$ and using Jessen's inequality, we have
\begin{eqnarray*}
\frac{\sum_{i=1}^{N}\mathbb{E}d(Y^{i},X^{i})}{N}&\leq&\sqrt
{\frac{2C_{1}\sum_{i=1}^{N}\mathbb{E}H_{i}(Y^{[1,i-1]})}{N}}+
\frac{\sum_{i=1}^{N}\sum_{j=1}^{i-1}c_{ij}\mathbb
{E}d(Y^{j},X^{j})}{N}
\\
&=& \sqrt{\frac{2C_{1}H(Q|P)}{N}}+\frac{\sum_{j=1}^{N-1}\mathbb
{E}d(Y^{j},X^{j})\sum_{i=j+1}^{N}c_{ij}}{N}
\\
&\leq&\sqrt{\frac{2C_{1}H(Q|P)}{N}}+\frac{r_1
\sum_{j=1}^{N}\mathbb{E}d(Y^{j},X^{j})}{N}
\end{eqnarray*}
the above inequality gives us
\[
W_{1,d_{L^{1}}}(Q,P)\leq\sum_{i=1}^{N}
\mathbb{E}d\bigl(Y^{i},X^{i}\bigr)\leq\sqrt{2
\frac{NC_{1}}{(1-r_{1})^{2}}H(Q|P)},
\]
that is,
$P=P(x_0,\cdot)\in T_{1} (\frac{NC_{1}}{(1-r_{1})^{2}} )$.
\end{pf}

Theorem \ref{concentrationinequalities} is based on the following
dependent tensorization result of Djellout, Guillin and Wu \cite{DGW}.

%le3.5 #&#
\begin{lem}[(\cite{DGW}, Theorem 2.11)]\label{lemDGW}  Let $\mathbb{P}$ be a probability measure on the
product space $(E^{n}, \mathscr{B}^{n}), n\geq2$. For any
$x=(x_1,\ldots,x_n)\in
E^{n}, x_{[1,k]}:=(x_{1},\ldots, x_{k})$. Let
$\mathbb{P}_{k}(\cdot|x_{[1,k-1]})$ denote the regular conditional
law of $x_{k}$ given $x_{[1,k-1]}$ under $\mathbb{P}$ for $2\le k\le
n$, and $\mathbb{P}_{k}(\cdot|x_{[1,k-1]})$
be the distribution of $x_1$ for $k=1$.

Assume that:
\begin{enumerate}[(2)]
\item[(1)] For some metric $d$ on $E$, $\mathbb{P}_{k}(\cdot
|x_{[1,k-1]})\in T_{1}(C)$ on $(E,
d)$ for all $k\geq1, x_{[1,k-1]}\in E^{k-1}$;

\item[(2)] there is some constant $S>0$ such that for all real bounded
Lipschitzian function $f(x_{k+1},\ldots,x_{n})$ with
$\|f\|_{\mathrm{Lip}(d_{L^{1}})}\leq1$, for all $x\in E^{n}, y_{k}\in E$,
\begin{eqnarray*}
&&\bigl|\mathbb{E}_{\mathbb{P}}\bigl(f(X_{k+1},\ldots,X_{n})|X_{[1,k]}=x_{[1,k]}
\bigr)-\mathbb{E}_{\mathbb{P}}\bigl(f(X_{k+1},\ldots,X_{n})
|X_{[1,k]}=(x_{[1,k-1]}, y_{k})\bigr)\bigr|\\
&&\quad
\leq Sd(x_{k},y_{k}).
\end{eqnarray*}
\end{enumerate}
Then for all function $F$ on $E^{n}$ satisfying $\|F\|
_{\mathrm{Lip}(d_{L^1})}\le\alpha$, we have
\[
\mathbb{E}_{\mathbb{P}}\RMe^{\lambda(F-\mathbb{E}_{\mathbb{P}}F)}\leq\exp
\bigl(C\lambda^{2}(1+S)^{2}
\alpha^{2}n/2\bigr)\qquad \forall\lambda\in\mathbb{R}. %
\]
Equivalently, $\mathbb{P}\in T_{1}(C_{n})$ on $(E^{n},
d_{L^{1}})$ with
\[
C_{n}=nC(1+S)^{2}.
\]
\end{lem}

We are now ready to prove Theorem \ref{concentrationinequalities}.

\begin{pf*}{Proof of Theorem \ref{concentrationinequalities}}
We
will apply Lemma \ref{lemDGW} with $(E,d)$ being $(E^N, d_{L^1})$, and
$\pp$ be the law of $(Z_1,\ldots,Z_n)$
on $(E^N)^n$.

By (\ref{matrixineq}), Lemma \ref{1-norm} and the condition that
$r_1<1/2$, the constant $S$ in Lemma \ref{lemDGW} is bounded from above
by
\[
\sup_{x,y\in
E^{N}}\frac{1}{d_{L^{1}}(x,y)}\mathbb{E}^{x,y}\sum
_{k=1}^{\infty
}d_{L^{1}}(X_{kN},Y_{kN})
\leq\sum_{k=1}^{\infty}\biggl(
\frac{r_{1}}{1-r_{1}}\biggr)^{k}=\frac{r_{1}}{1-2r_{1}}.
\]

Take $S=\frac{r_{1}}{1-2r_{1}},
F(Z_{1},\ldots,Z_{n})=\frac{1}{n}\sum_{k=1}^{n}f(Z_{k})$, then the
Lipschitzian norm $\|F\|_{\mathrm{Lip}}$ of $F$ w.r.t. the
$d_{L^1}(x,y)=\sum_{k=1}^n d_{L^1}(x_k,y_k)$ (for $x,y\in(E^N)^n$)
is not greater than $\|f\|_{\mathrm{Lip}(d_{L^{1}})}/n\le\alpha/n$. Thus by
Lemmas \ref{lemDGW} and \ref{transine},
\begin{eqnarray*}
&&
\mathbb{E}^x \exp\Biggl(\lambda\Biggl(\frac{1}{n}\sum
_{k=1}^{n}f(Z_{k})-
\frac
{1}{n}\sum_{k=1}^{n}P^{k}f(x)
\Biggr) \Biggr)\\
&&\quad\leq\exp\biggl\{\frac{NC_{1}}{(1-r_{1})^{2}}\lambda
^{2}(1+S)^{2}
\biggl(\frac
{\alpha}{n}\biggr)^{2}n/2 \biggr\}\qquad \forall\lambda\in
\mathbb{R}.
\end{eqnarray*}
So, by the classic approach, firstly using
Chebyshev's inequality, and then optimizing over $\lambda\geq0$, we
obtain the desired part $(a)$ in Theorem \ref{concentrationinequalities}.

Furthermore by Theorem \ref{convergencerate}, we have
\begin{eqnarray*}
\Biggl|\frac{1}{n}\sum_{k=1}^{n}P^{k}f(x)-
\mu(f)\Biggr|&\leq&\frac{1}{n}\sum_{k=1}^{n}\bigl|P^{k}f(x)-
\mu(f)\bigr|
\\
&\leq&\frac{1}{n}\sum_{k=1}^{n}
r^{k}\max_{1\leq i\leq
N}\int_{E^{N}}d
\bigl(x^{i},y^{i}\bigr)\mrmd \mu(y)\cdot\sum
_{i=1}^N \delta_i(f)
\\
&\leq&\frac{1}{n}\frac{ r}{1-r}\max_{1\leq i\leq
N}\int
_{E^{N}}d\bigl(x^{i},y^{i}\bigr)\mrmd \mu(y)
\cdot\sum_{i=1}^N \delta_i(f).
\end{eqnarray*}
Thus, we obtain part $(b)$ in Theorem \ref{concentrationinequalities}
from its part (a).
\end{pf*}

% zodis "Acknowledgments" paliekamas pagal autoriu
\section*{Acknowledgements}

Supported in part by Thousand Talents Program of the Chinese Academy
of Sciences and le projet ANR EVOL.
We are grateful to the two referees for
their suggestions and references, which improve sensitively the
presentation of the paper.

%suskaldyti doi

% imsref loaded by lrinkeviciute, 2013-08-16 08:36:47
% imsref loaded by lrinkeviciute, 2013-08-20 13:19:40

\printhistory

\end{document}